\documentclass[twocolumn,amsmath,amssymb,aps,secnumarabic,%
    nofootinbib,groupedaddress]{revtex4}

\usepackage{graphics,times,mathptmx,amsthm,pstricks}
\usepackage[small]{subfigure}
\newtheorem{theorem}{Theorem}
\newtheorem{conjecture}[theorem]{Conjecture}
\newtheorem{lemma}[theorem]{Lemma}

\newcommand{\F}{\mathbb{F}}
\newcommand{\hF}{\widehat{F}}
\newcommand{\hh}{\widehat{h}}
\newcommand{\tX}{\widetilde{X}}
\newcommand{\defeq}{\stackrel{\mathrm{def}}{=}}
\newcommand{\deq}{\stackrel{\mathrm{d}}{=}}

\newcommand{\Eg}{\emph{E.g.}}
\newcommand{\ie}{\emph{i.e.}}

\newcommand{\thm}[1]{Theorem~\ref{#1}}
\newcommand{\lem}[1]{Lemma~\ref{#1}}
\newcommand{\Step}[1]{\noindent \emph{Step #1.}}

\renewcommand{\sec}[1]{Section~\ref{#1}}
\newcommand{\fig}[1]{Figure~\ref{#1}}
\newcommand{\eq}[2]{\begin{equation}\label{#1}#2\end{equation}}
\newcommand{\Eqr}[1]{Equation~\eqref{#1}}
\newcommand{\eqr}[1]{equation~\eqref{#1}}
\newcommand{\eqpr}[1]{eq.~\eqref{#1}}
\newcommand{\eatline}{\vspace{-\baselineskip}}
\newcommand{\eathalfline}{\vspace{-.5\baselineskip}}

\newenvironment{fullfigure}[2]
    {\begin{figure}[htb]\begin{center}\def\fullfiga{#1}\def\fullfigb{#2}}
    {\vspace{\baselineskip}\caption{\fullfigb.}\label{\fullfiga}
        \end{center}\end{figure}}
\newenvironment{fullfigure*}[2]
    {\begin{figure*}[htb]\begin{center}\def\fullfiga{#1}\def\fullfigb{#2}}
    {\vspace{\baselineskip}\caption{\fullfigb.}\label{\fullfiga}
        \end{center}\end{figure*}}

\psset{linewidth=.5pt,dash=3pt 3pt,doublesep=.05,dotsize=1pt 5,arrowsize=2pt 3}
\newcommand{\ptex}{\psline[linewidth=1pt](.1,.1)(-.1,-.1)
    \psline[linewidth=1pt](.1,-.1)(-.1,.1)}

\begin{document}
\title{Special moments}

\author{Greg Kuperberg}
\email{greg@math.ucdavis.edu}
\thanks{Supported by NSF grant DMS \#0306681}
\affiliation{Department of Mathematics, University of
    California, Davis, CA 95616}

\begin{abstract}
\centerline{\textit{\normalsize Dedicated to the memory of David Robbins}}
\vspace{\baselineskip}

In this article, we show that a linear combination $\tX$ of $n$ independent,
unbiased Bernoulli random variables $\{X_k\}$ can match the first $2n$ moments
of a random variable $Y$ which is uniform on an interval. More generally, for
each $p \ge 2$, each $X_k$ can be uniform on an arithmetic progression of
length $p$.  All values of $\tX$ lie in the range of $Y$, and their ordering as
real numbers coincides with dictionary order on the vector $(X_1,\ldots,X_n)$.

The construction involves the roots of truncated $q$-exponential series. It
applies to a construction in numerical cubature using error-correcting codes
\cite{Kuperberg:cubature}. For example, when $n=2$ and $p=2$, the values of
$\tX$ are the 4-point Chebyshev quadrature formula.
\end{abstract}

\maketitle

\section{Introduction}
\label{s:intro}

One of the standard proofs of the central limit theorem establishes that the
moments of the normalized sum
$$\frac{X_1 + X_2 + \ldots + X_n}{\sqrt{n}}$$
of $n$ i.i.d. centered random variables with finite moments converge to the
moments of a Gaussian random variable.  This fact raises the question of when
the first $n$ moments of a random variable $Y$ can be matched by a linear
combination  of independent copies of another variable $X$. It is an easy
exercise with cumulants that this is impossible when $Y$ is Gaussian and $X$ is
not Gaussian.  In this article we will show that for every $n$, $Y$ can be the
uniform distribution on an interval if $X$ is unbiased Bernoulli.  More
generally $X$ can be uniform on an arithmetic progression of length $p$ for any
$p \ge 2$.

We conjecture that the moments of most absolutely continuous distributions
cannot be matched by those of a linear combination of Bernoulli random
variables.  In this sense the uniform distribution on an interval has ``special
moments''.

\begin{theorem} Let $p \ge 2$ be an integer, let $X$ be a uniformly random
variable on the set
$$\{p-1,p-3,p-5,\ldots,1-p\},$$
and let  $X_1,X_2,\ldots,X_n$ be independent copies of $X$. Then there exist
unique constants
$$a_1 > a_2 > \cdots > a_n > 0$$
such that the first $2n$ moments of
$$\tX = a_1 X_1 + a_2 X_2 + \cdots + a_n X_n$$
agree with the first $2n$ moments of a random variable $Y$ which is uniform on
$[-1,1]$. Moreover
$$\sum_{j=1}^n |a_j-p^{-j}| < \frac{1}{p^n(p-1)}.$$
\label{th:main} \end{theorem} \eatline

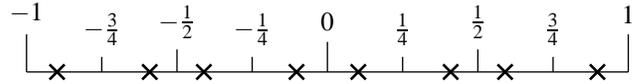
\begin{fullfigure}{f:ruler}{The values of $\tX$ marked on a ruler
    when $p=2$ and $n=3$}
\pspicture(-4,0)(4,1.25)
\psline(-4,0)(4,0)\psline(-4,0)(-4,.5)\psline(4,0)(4,.5)\psline(0,0)(0,.4)
\psline(-2,0)(-2,.3)\psline(2,0)(2,.3)
\psline(-3,0)(-3,.2)\psline(-1,0)(-1,.2)\psline(1,0)(1,.2)\psline(3,0)(3,.2)
\rput(3.592044,0){\ptex} \rput(-3.592044,0){\ptex}
\rput(2.360508,0){\ptex} \rput(-2.360508,0){\ptex}
\rput(1.640516,0){\ptex} \rput(-1.640516,0){\ptex}
\rput(0.408980,0){\ptex} \rput(-0.408980,0){\ptex}
\rput[b](-4,.65){\normalsize $-1$}
\rput[b](-3,.35){\normalsize $-\frac34$}
\rput[b](-2,.45){\normalsize $-\frac12$}
\rput[b](-1,.35){\normalsize $-\frac14$}
\rput[b](0,.55){\normalsize $0$}
\rput[b](1,.35){\normalsize $\frac14$}
\rput[b](2,.45){\normalsize $\frac12$}
\rput[b](3,.35){\normalsize $\frac34$}
\rput[b](4,.65){\normalsize $1$}
\endpspicture\eathalfline
\end{fullfigure}

The second part of \thm{th:main} implies that all $p^n$ values of $\tX$ lie in
the interior of the interval $[-1,1]$.  Moreover, if we cut $[-1,1]$ into
$p^n$ equal subintervals, then each subinterval contains exactly one value of
$\tX$.  This pattern is
summarized by the relation
$$\tX \in \bigl(\sum_{j=1}^n X_jp^{-j}-p^{-n},
    \sum_{j=1}^n X_jp^{-j}+p^{-n}\bigr) \subseteq (-1,1),$$
which is easily equivalent to the inequality in \thm{th:main}.

A more combinatorial formulation is as follows:  If we number equal
subintervals of $[-1,1]$ from $0$ to $p^n-1$ and if a particular value of $\tX$
lies in subinterval $k$, then the $j$th digit of $k$ in base $p$ is determined
by the $j$th variable $X_j$.  In symbols,
$$k=\sum_{j=1}^n \frac{X_j+p-1}{2} p^{n-j}.$$

\fig{f:ruler} shows an example.  If $p=2$ and $n=3$, then
$$\tX \approx \pm .500128 \pm .243941 \pm .153942,$$
where the signs are independent and random.   It is natural to think of
$[-1,1]$ with its $p^n$ subintervals as a $p$-adic ruler. In \fig{f:ruler}, the
8 values of $\tX$ with $p=2$ and $n=3$ are marked on a dyadic ruler.  The
figure bears out \thm{th:main} because each mark lies in a different gap
between ruler marks. The marks appear to be at a fixed distance from the ruler
marks, but this is only approximately true:  The first two terms of $\tX$ are
close to $\pm \frac12$ and $\pm \frac14$, but not exactly equal.  In
the general case, \lem{l:bound} below establishes that the coefficient $a_j$ is
extremely close to the (halved) ruler spacing $p^{-j}$ when $j$ is small:
\eq{e:k2}{a_jp^j = 1 \pm O(p^{-k^2-k}),}
where $k = n+1-j$.

Let $Z_{p,n}$ be the range of $\tX$.  Since $\tX$ has the same first $2n$
moments as $Y$, indeed trivially the same $(2n+1)$st moment as well,
the equation
$$\frac12\int_{-1}^1 P(x) dx = \frac1{p^n} \sum_{x \in Z_{p,n}} P(x)$$
holds for any polynomial $P$ of degree at most $2n+1$. A weighted set $Z$ with
this property up to some degree $t$ is called an (interpolatory)
\emph{$t$-quadrature formula}.  \Eg, Simpson's rule and Gaussian quadrature are
standard quadrature formulas.  Our quadrature formula $Z_{p,n}$ is highly
inefficient for general $p$ and $n$, but its special structure is useful for
the higher-dimensional cubature problem for integration on the $k$-cube
$[-1,1]^k$. Elsewhere \cite{Kuperberg:cubature} we combine the product formula
$Z_{2,n}^{\times d}$ with binary error-correcting codes, in particular extended
BCH codes, to obtain a $(2n+1)$-cubature formula on $[-1,1]^k$ with equal
weights and $O(k^n)$ points.  (The asymptotic bound is with $n$ fixed and $k
\to \infty$.)

\section{The proof}
\label{s:proof}

We will write $a_{n,j}$ for $a_j$, to make clear that they depend on $n$.

\begin{lemma} The random variables $\tX$ and $Y$ have the same first $2n$
moments if and only if
$$\sum_{j=1}^n a_{n,j}^{2k} = \frac1{p^{2k}-1}$$
for all $1 \le k \le n$.
\label{l:powsum} \end{lemma}

\begin{proof}
Recall that $Y$ is uniformly random on $[-1,1]$ and that $X$ is uniformly
random on
$$\{p-1,p-3,p-5,\ldots,1-p\}.$$
Thus if they are independent, then $X+Y$ is uniformly random on $[-p,p]$, so
\eq{e:xy}{X+Y \deq pY.}
Here ``$\deq$'' means equality of distribution.

To understand the implications of this relation between $X$ and $Y$, we review
the relations among cumulants, moments, and their generating functions
\cite{AS:handbook}.  For a general random variable $X$, the $k$th moment is
denoted $\mu_k(X)$, the exponential generating function of all of the moments
is the moment function $M_X(t)$, the cumulant function $K_X(t)$ is its
logarithm, and the cumulants $\kappa_k(X)$ are defined by $K_X(t)$ as their
exponential generating function.  In formulas,
\begin{align*}
\mu_k(X) &\defeq E[X^k]
& M_X(t) \defeq E[e^{Xt}] &= \sum_{k=0}^\infty \frac{\mu_k(X)t^k}{k!} \\
K_X(t) &\defeq \ln M_X(t) & \sum_{k=0}^\infty \frac{\kappa_k(X)t^k}{k!}
    &\defeq K_X(t)
\end{align*}
This framework is designed so that first, cumulants carry the same information
as moments, and second, cumulants are additive, \ie,
$$K_X(t) + K_Y(t) = K_{X+Y}(t),$$
for independent random variables $X$ and $Y$.

\Eqr{e:xy} yields the cumulant generating function equation
$$K_X(t) + K_Y(t) = K_Y(pt),$$
which we can write as a relation between individual cumulants:
$$\kappa_k(X) + \kappa_k(Y) = p^k\kappa_k(Y).$$
The odd cumulants of $X$ and $Y$ vanish since they are even random variables.
The even cumulants thus satisfy the relation
$$\kappa_{2k}(X) = (p^{2k}-1)\kappa_{2k}(Y).$$
Since
$$\kappa_{2k}(\tX) = \sum_j a_{n,j}^{2k}\kappa_{2k}(X),$$
it suffices for the $a_{n,j}$'s to satisfy the stated power sum relation.

This condition is also necessary provided that each $\kappa_{2k}(Y) \ne 0$. To
check this, we will establish that $\kappa_{2k}(X) \ne 0$ when $p=2$. The
moment function with imaginary argument, $M_X(it)$, is also called the
characteristic function of $X$ (meaning the Fourier transform of the
distribution of $X$).  In this case, its logarithm is:
$$K_X(it) = \log \cos t \qquad K_X(it)' = -\tan t.$$
The relation
$$(\tan t)' = (\tan t)^2 + 1$$
implies that the tangent function has strictly positive odd derivatives,
so
$$(-1)^{k+1}\kappa_{2k}(X) > 0$$
for all $k$.
\end{proof}

For convenience let $q = p^2$, and let:
$$b_{n,j} = a_{n,j}^2 \qquad r_{n,j} = \frac1{b_{n,j}}.$$
\lem{l:powsum} can then be restated as
$$\sum_{j=1}^n b_{n,j}^k = \frac1{q^k-1}$$
for all $1 \le k \le n$.  This implies a unique solution for the $a_{n,j}$'s
provided that each $b_{n,j}$ is real and positive.  For convenience we will
study a polynomial whose roots are $r_{n,j}$ for $1 \le j \le n$.

\begin{lemma} If
$$\sum_{j=1}^n b_{n,j}^k = \frac1{q^k-1}$$
for all $1 \le k \le n$, then
$$F_{q,\le n}(x) \defeq \prod_{j=1}^n (1-b_{n,j} x)
    = \sum_{k=0}^n \frac{x^k} {(1-q)(1-q^2)\cdots (1-q^k)}.$$
\eatline \label{l:poly} \end{lemma}

\begin{proof} Let
$$p_k = \sum_{j=1}^n b_{n,j}^k$$
be the $k$th power sum of the $b_{n,j}$'s, and let $e_k$ be the corresponding
elementary symmetric function, so that
$$\prod_{j=1}^n (1+b_{n,j} x) = 1+\sum_{k=1}^n e_k x^k.$$
Since the first $n$ elementary symmetric functions determine the first $n$
power sums and vice versa, and since our desired value
$$p_k = \frac1{q^k-1}$$
does not depend on $n$, we can derive each $e_k$ by taking the limit $n \to
\infty$ and finding $b_{\infty,j}$'s to match all $p_k$'s.  Let
$$b_{\infty,j} = q^{-j}.$$
Then
$$p_k = \sum_{j=1}^\infty b_{\infty,j}^k = \frac{1}{q^k-1}$$
since the left side is a geometric series.  Moreover
\begin{align*}
e_k &= \sum_{1 \le j_1 < j_2 < \cdots < j_k} \hspace{-1.25em}
    b_{\infty,j_1}b_{\infty,j_2} \cdots b_{\infty,j_k} \\
    &= \frac{1}{(q-1)(q^2-1)\cdots(q^k-1)}
\end{align*}
by a routine combinatorial exercise.  Another way to recognize these values of
$p_k$ and $e_k$ is that they are the principal specialization of the ring
$\Lambda$ of symmetric functions \cite[\S7.8]{Stanley:enumerative2},
transported by the fundamental involution $\omega$ and a sign involution
$\sigma$:
$$\omega(e_k) = h_k \qquad \sigma(e_k) = (-1)^k e_k.$$
To conclude, our explicit choice for the $b_{\infty,j}$'s establishes that the
given $p_k$'s are consistent with the claimed $e_k$'s.
\end{proof}

To continue the example mentioned first in \sec{s:intro},
$$F_{4,\le 3}(x) = 1-\frac{x}{3} + \frac{x^2}{45} - \frac{x^3}{2835}.$$
Its roots are
$$(\frac1{a_{3,1}^2}, \frac1{a_{3,2}^2}, \frac1{a_{3,3}^2}) \approx
    (3.997956, 16.80465, 42.19739),$$
so that
$$(a_{3,1},a_{3,2},a_{3,3}) \approx (.500128, .243941, .153942).$$
Therefore
$$\tX \approx \pm .500128 \pm .243941 \pm .153942$$
when $p=2$ and $n=3$.

We will need the $q$-Pochhammer symbol \cite{GR:hypergeometric}:
$$(a;q)_k = \begin{cases}
    \prod_{j=0}^{k-1} (1-aq^j) & k > 0 \\
    \prod_{j=1}^{-k} (1-aq^{-j})^{-1} & k < 0. \\
    1 & k = 0 \end{cases}$$

We also define
\begin{align}
F_{q,n}(x) &= \frac{x^n}{(q;q)_n} &
F_q(x) &= \sum_{k=0}^\infty \frac{x^k}{(q;q)_k} \nonumber \\
F_{q,\le n}(x) &= \sum_{k=0}^n \frac{x^k}{(q;q)_k} &
F_{q,>n}(x) &= \sum_{k={n+1}}^\infty \frac{x^k}{(q;q)_k}. \label{e:fdef}
\end{align}
(The polynomial $F_{q,\le n}(x)$ was already used in \lem{l:poly}.) We will
often use the relation
\eq{e:ratio}{\frac{F_{q,n}(x)}{F_{q,n-1}(x)} = \frac{x}{1-q^n}.}

The function $F_q(x)$ is related to the standard Jackson $q$-exponential
function $e_q(x)$ \cite{GR:hypergeometric} by
$$F_q(x) = e_q(\frac{x}{1-q}).$$
(In some works, $F_q(x)$ itself is called a $q$-exponential.)
By the proof of \lem{l:poly},
\eq{e:fprod}{F_q(x) = \prod_{j=1}^\infty (1-q^{-j}x)
    = \frac{1}{(x;q)_{-\infty}}.}
This identity holds for all $q$ as an equality of formal power series. We will
need the stronger fact that it is an equality of entire analytic functions
when $q > 1$.

To establish \thm{th:main}, we would like to understand the effect of
truncation on the first $n$ zeroes of $F_q(x)$ when $q \ge 4$.

\begin{lemma} If $q \ge 4$, then $F_{q,\le n}(x)$ has $n$ distinct, positive
roots.
\label{l:unimodal}
\end{lemma}

\begin{proof} It suffices to show that the value
$$f_{k,\le n} = F_{q,\le n}(q^{k+\frac12})$$
alternates in sign as $k$ ranges from $0$ to $n$.  The terms of $f_{k,\le n}$
are
$$f_{k,j} = \frac{q^{j(k+\frac12)}}{(q;q)_j}.$$
These alternate in sign in $j$ and we claim that $f_{k,\le n}$ has the same
sign as $f_{k,k}$.  This claim will imply the lemma.

By \eqr{e:ratio}, the sequence $\{|f_{k,j}|\}$ is unimodal in $j$ and achieves
its maximum at $j=k$.  This already implies that $f_{k,\le n}$ has the same
sign as $f_{k,k}$ if $k=0$ or $k=n$.  If $1 < k < n$, then
$$\bigl|\sum_{j=0}^{k-1} f_{k,j}\bigr| < |f_{k,k-1}| \qquad
    \bigl|\sum_{j=k+1}^n f_{k,j}\bigr| < |f_{k,k+1}|.$$
Finally
$$\frac{|f_{k,k-1}| + |f_{k,k+1}|}{|f_{k,k}|} = \frac{q^k-1}{q^{k+\frac12}}
     + \frac{q^{k+\frac12}}{q^{k+1}-1} < \frac{2}{q^{\frac12}} \le 1.$$
The first equality once again comes from \eqr{e:ratio}, while the
last inequality is the only step that requires $q \ge 4$ instead of merely
$q > 1$.  Thus in each case $f_{k,k}$ is the dominant term in $f_{k,\le n}$.
\end{proof}

We continue our example case with $q=4$ and $n=3$ to illustrate
\lem{l:unimodal} and its proof:
\begin{align*}
F_{4,\le 3}(2) &= 1-\frac{2}3 + \frac{4}{45} - \frac{8}{2835}
    = \frac{1189}{2835} \\
F_{4,\le 3}(8) &= 1-\frac{8}3 + \frac{64}{45} - \frac{512}{2835}
    = -\frac{1205}{2835} \\
F_{4,\le 3}(32) &= 1-\frac{32}3 + \frac{1024}{45} - \frac{32768}{2835}
    = \frac{4339}{2835} \\
F_{4,\le 3}(128) &= 1-\frac{128}3 + \frac{16384}{45} - \frac{2097152}{2835}
    = -\frac{1183085}{2835}.
\end{align*}
In this example we can see that the $k$th sum is dominated by its $k$th term.

  Since
$F_{4,\le 3}(x)$ is a cubic polynomial with four values that alternate in sign,
it therefore has distinct, positive roots, previously noted to be
$$(r_{3,1}, r_{3,2}, r_{3,3}) \approx (3.997956, 16.80465, 42.19739).$$

\lem{l:unimodal} is also implied by the final lemma, \lem{l:bound}, but the
proof of \lem{l:bound} is much more complicated.  It may unfortunately be as
taxing for the reader as it was for the author.

\begin{lemma} Let $q \ge 4$.  Let $r_{n,j}$
$$r_{n,1} < r_{n,2} < \cdots < r_{n,n}$$
be the roots of $F_{q,\le n}(x)$
and let $k = n+1-j$.  Then $r_{n,j}q^{-j}$ lies between $1$ and
\eq{e:ck}{c_k = \begin{cases} 1-2q^{-1} & k=1 \\
    1+(-1)^k4q^{-\binom{k+1}2} & k>1 \end{cases}}
for every $1 \le j \le n$.
\label{l:bound} \end{lemma}

\begin{proof} Let
$$\hF_{q,m}(x) = (-1)^{n+1}F_{q,m}(x)$$
with $m$ independent of $n$, and extend this notation to all of the
definitions in \eqref{e:fdef}.

The proof consists of three steps.  In the first step, we will show that
the lemma follows from the inequality
\eq{e:fhat}{\hF_q(c_kq^j) > \hF_{q,>n}(c_kq^j)}
and we will show that both sides are positive.  In the second step,
we will show that this inequality follows from the inequality
$$h_{q,k} \defeq \lim_{n \to \infty}
    \frac{\hF_q(c_kq^j)}{\hF_{q,>n}(c_kq^j)} \ge 1,$$
where $k$ is fixed in taking the limit.  Finally in the third step, we
will show that
$$h_{q,k} > 1.$$

\Step{1}  We claim that \eqr{e:fhat} implies that \eq{e:fdiff}{\hF_{q,\le n}(x)
= \hF_q(x) - \hF_{q,>n}(x)} changes sign as $x$ passes from $q^j$ to $c_kq^j$. 
In the estimates for this step we will assume that
$$q^j \le x \le c_kq^j \quad \mbox{or} \quad c_kq^j \le x \le q^j,$$
except when we explicitly state more general conditions.

By \eqr{e:fprod},
$\hF_q(q^j)$ vanishes when $j \ge 1$.  
When $x>0$,
$$\hF_{q,n+1}(x) > 0$$
by the definition of $\hF_{q,n+1}(x)$.  Moreover
\eq{e:pos}{\hF_{q,>n}(x) > 0}
because, by \eqr{e:ratio}, the series for $\hF_{q,>n}(x)$ is alternating and
decreasing when $0 < x < q^{n+2}-1$.  Thus
$$\hF_{q,>n}(q^j) > 0$$
since $j \le n$ and the argument $x = q^j$ is thus in the required range.
\Eqr{e:fdiff} now tells us both that
$$\hF_{q,\le n}(q^j) < 0$$
willy-nilly, and that
$$\hF_{q,\le n}(c_kq^j) > 0$$
is equivalent to \eqr{e:fhat}.  This establishes the first claim of this step.

\Eqr{e:pos} also tells us that
$$\hF_{q,>n}(c_kq^j) > 0,$$
since $x = c_kq^j$ is also in the required range $0 < x < q^{n+2}-1$. Finally
we confirm that
$$\hF_q(c_kq^j) > 0$$
by counting sign changes in the factors of \eqr{e:fprod}.  To do this properly,
observe from \eqr{e:ck} that $q^{j-1} < c_kq^j < q^j$ when $k$ is odd and $q^j
< c_kq^j < q^{j+1}$ when $k$ is even.  Either way, the number of sign changes
has the same parity as $j+k = n+1$.

\Step{2} The goal of this messy step is to reduce \eqr{e:fhat} to its
asymptotic limit $n,j \to \infty$ with $k$ fixed.  We will use some preliminary
relations for the Pochhammer symbol.  The  product relation
\eq{e:pochprod}{(a;q)_m = (a;q)_\ell(aq^\ell;q)_{m-\ell}}
holds when $\ell$ is finite (but $m$ need not be).  The
inversion relation
\eq{e:pochinv}{(aq;q)_\ell
    = \frac{(-a)^\ell q^{\binom{\ell+1}2}}{(a^{-1};q)_{-\ell}}}
holds for all finite $\ell$; it follows from the trivial identity
$$1-a = -a(1-a^{-1}).$$
The inequality
\eq{e:pochle}{(a;q)_{-\ell} \ge 1}
holds when $0 \le \ell \le \infty$ and $0 < a < q$ (with equality only when
$\ell = 0$).  Finally
\eq{e:pochdiff}{(a;q)_{\ell+1} - (a;q)_\ell = -aq^\ell(a;q)_\ell}
for all finite $\ell$. We will also use the elementary binomial identity
\eq{e:eb}{\binom{\ell+m}{2} = \binom{\ell}2 + \ell m + \binom{m}2.}

The left side of \eqr{e:fhat} limits to a product of manageable factors:
\begin{align*}
\hF_q(c_kq^j)
    &= \frac{(-1)^{n+1}}{(c_kq^j;q)_{-\infty}} \\
    &= \frac{(-1)^{n+1}(1-c_k)(c_kq;q)_{j-1}}{(c_k;q)_{-\infty}}
        & \mbox{(by \eqpr{e:pochprod})} \\
    &= \frac{(-1)^{k+1}(1-c_k)q^{\binom{j}2}c_k^{j-1}}
        {(c_k^{-1};q)_{1-j}(c_k;q)_{-\infty}}
        & \mbox{(by \eqpr{e:pochinv})} \\
    &\ge \frac{(-1)^{k+1}(1-c_k)q^{\binom{j}2}c_k^{j-1}}
        {(c_k^{-1};q)_{-\infty}(c_k;q)_{-\infty}}
        & \mbox{(by eqs. (\ref{e:pochprod},\ref{e:pochle})).}
\end{align*}

Meanwhile, the right side of \eqr{e:fhat} essentially stabilizes
as a power series in $xq^{-n-2}$:
\begin{align*}
\hF_{q,>n}(x)
    &= (-1)^{n+1} \sum_{\ell = n+1}^\infty \frac{x^\ell}{(q;q)_\ell} \\
    &= (-x)^{n+1} \sum_{\ell = 0}^\infty \frac{x^\ell}{(q;q)_{\ell+n+1}} \\
    &= x^{n+1} \sum_{\ell = 0}^\infty \frac{(-x)^\ell
        (1;q)_{-\ell-n-1}}{q^{\binom{\ell+n+2}{2}}}
        & \hspace{-.5em} \mbox{(by \eqpr{e:pochinv})} \\
    &= \frac{x^{n+1}}{q^{\binom{n+2}{2}}} \sum_{\ell = 0}^\infty
        \frac{(-xq^{-n-2})^\ell (1;q)_{-\ell-n-1}}{q^{\binom{\ell}{2}}}.
        & \hspace{-.5em} \mbox{(by \eqpr{e:eb})} \\
\end{align*}
To isolate the power series, let
$$G_{q,n}(t) = \sum_{\ell=0}^{\infty}
     \frac{(-t)^\ell(1;q)_{-\ell-n-1}}{q^{\binom{\ell}{2}}}.$$
Then
$$\hF_{q,>n}(x) = \frac{x^{n+1}}{q^{\binom{n+2}2}}G_{q,n}(xq^{-n-2}).$$

To complete step 2, we will show that $G_{q,n}(x)$ is
monotonic in $n$ and consolidate inequalities.  Observe that
$$G_{q,m+1}(t) - G_{q,m}(t) = \sum_{\ell=0}^{\infty}
     \frac{(-t)^\ell (1;q)_{-\ell-m-2}}{q^{\binom{\ell}{2}+\ell+m+2}}$$
by \eqr{e:pochdiff}.  This series is alternating decreasing when
$0<t<q-q^{-m-2}$ and $m>0$, whence
$$G_{q,m+1}(t) > G_{q,m}(t).$$
In particular,
$$G_{q,\infty}(t) \defeq \sum_{\ell=0}^{\infty}
    \frac{(-t)^\ell(1;q)_{-\infty}}{q^{\binom{\ell}{2}}} > G_{q,n}(t)$$
when $0<t<q-q^{-n-2}$, by induction on $m \ge n$. In particular
$$\hF_{q,>n}(c_kq^j) < \frac{(c_kq^j)^{n+1}}{q^{\binom{n+2}2}}
    G_{q,\infty}(c_kq^{-k-1})$$
because $t = c_kq^{-k-1}$ is well within range given that $k \ge 1$,
$c_k < q$, and $q \ge 4$.

Thus we have estimates for both sides of \eqr{e:fhat}.  Upon
close examination, they sacrifice less and less as $n \to \infty$.  Combining
the estimates,
\eq{e:hqk}{\frac{\hF_q(c_kq^j)}{\hF_{q,>n}(c_kq^j)} \stackrel{>}{\to}
    \frac{(-1)^{k+1}(1-c_k)q^{\binom{k+1}2}c_k^{-k-1}}
    {(c_k^{-1};q)_{-\infty}(c_k;q)_{-\infty}G_{q,\infty}(c_kq^{-k-1})}.}
The right side is the limit $h_{q,k}$ defined previously.  We have
shown that the lemma follows from the inequality $h_{q,k} \ge 1$.
It is also necessary if our use of the intermediate
value theorem in step 1 is to work for all $n$.

\Step{3} We will need the Pochhammer symbol estimate
\eq{e:pochgeom}{(a;q)_{-\infty}^{-1} > 1-\frac{a}{q-1}.}
We claim that it holds when $0 < a < q$ and $q > 1$.  By \eqr{e:pochle}
$$(a;q)_{-\ell}^{-1} = (1 - aq^{-\ell})(a;q)^{-1}_{1-\ell}
    \ge (a;q)^{-1}_{1-\ell} - aq^{-\ell}$$
for any $\ell \ge 1$, with equality only when $\ell = 1$.  Thus by induction,
$$(a;q)_{-\ell}^{-1} \ge 1 - \sum_{m=1}^\ell aq^{-m},$$
again with equality only when $\ell = 1$.  Now we sum the geometric series in
the limit $\ell \to \infty$.

We will also need the estimate
\eq{e:gest}{G_{q,\infty}(t) < (1;q)_{-\infty},}
which holds when $0 < t < 1$ because the power series for $G_{q,\infty}(t)$ is
then alternating decreasing.  We apply equations \eqref{e:pochgeom} and
\eqref{e:gest} to the right side of \eqr{e:hqk} to obtain
\begin{align*}
h_{q,k} > &(-1)^{k+1}(1-c_k)q^{\binom{k+1}2}c_k^{-1-k} \nonumber \\
    &\cdot \bigl(1-\frac{c_k}{q-1}\bigr)\bigl(1-\frac{c_k^{-1}}{q-1}\bigr)
    \bigl(1-\frac{1}{q-1}\bigr) \defeq \hh_{q,k}
\end{align*}
when $k \ne 2$.  The inequality holds when $k=2$ as well, but it
is not adequate because our proof is a close shave in this case.
So we will define $\hh_{q,2}$ differently.  We refine \eqr{e:pochgeom},
$$(a;q)_{-\infty}^{-1} = \bigl(1-\frac{a}{q}\bigr)(aq^{-1};q)_{-\infty}^{-1}
    > \bigl(1-\frac{a}{q}\bigr) \bigl(1-\frac{a}{q(q-1)}\bigr),$$
to obtain
\begin{align*}
h_{q,2} > &(c_2-1)q^3 c_2^{-3} \bigl(1-\frac{c_2}{q-1}\bigr)
    \bigl(1-\frac{c_2^{-1}}{q-1}\bigr) \nonumber \\
    & \qquad \cdot \bigl(1-\frac{1}{q}\bigr)\bigl(1-\frac{1}{q(q-1)}\bigr)
    \defeq \hh_{q,2}.
\end{align*}

If we apply \eqr{e:ck} to the first occurrence of $c_k$ here, we learn that
$$\hh_{q,k} = 4c_k^{-1-k}\bigl(1-\frac{c_k}{q-1}\bigr)
    \bigl(1-\frac{c_k^{-1}}{q-1}\bigr) \bigl(1-\frac{1}{q-1}\bigr)$$
when $k>2$, while
\begin{align*}
\hh_{q,1} &= 2c_1^{-2}\bigl(1-\frac{c_1}{q-1}\bigr)
    \bigl(1-\frac{c_1^{-1}}{q-1}\bigr) \bigl(1-\frac{1}{q-1}\bigr) \\
\hh_{q,2} &= 4c_2^{-3}\bigl(1-\frac{c_2}{q-1}\bigr)
    \bigl(1-\frac{c_2^{-1}}{q-1}\bigr) \bigl(1-\frac{1}{q}\bigr)
    \bigl(1-\frac{1}{q(q-1)}\bigr)
\end{align*}

We claim that $\hh_{q,k} > 1$.  It can be checked directly with symbolic
algebra that
$$\hh_{q,1} = \frac{2q(q^2-4q+2)(q^2-2q+2)}{(q-1)^3(q-2)^2} > 1$$
when $q \ge 4$ (indeed when $q \ge 3.718$), and that
$$\hh_{q,2} = \frac{4q^4(q^2\!-\!q\!-\!1)(q^2\!-\!2q\!+\!2)
    (q^2\!-\!2)(q^4\!-\!2q^3\!-\!4)}{(q-1)^2(q^3+4)^4} > 1$$
when $q \ge 4$ (indeed when $q \ge 3.974$). When $k>2$, we claim that
\begin{align*}
c_k^{-1-k} &> \frac{99}{100} & 1-\frac{1}{q-1} &\ge \frac23 \\
1-\frac{c_k}{q-1} &> \frac58 & 1-\frac{c_k^{-1}}{q-1} &> \frac58,
\end{align*}
To check the first of these inequalities,  we apply \eqr{e:ck} and take the
logarithm of both sides.  We want to show that
$$(1+k)\log \bigl(1+(-1)^k4q^{-\binom{k+1}2}\bigr) < \log \frac{100}{99}.$$
We can assume that $k$ is even, so that $k \ge 4$.  We can
simplify using the elementary inequalities
$$\log (1+x) < x \qquad \log \frac{1}{1-x} > x.$$
Thus it suffices to show that
$$(1+k)4q^{-\binom{k+1}2} < \frac{1}{100}.$$
This holds easily assuming that $q \ge 4$ and $k \ge 4$.  The other three
inequalities also hold easily given that $q \ge 4$.

Thus when $k>2$,
$$\hh_{q,k} \ge 4 \cdot \frac{99}{100}
    \cdot \frac58 \cdot \frac58 \cdot \frac23 = \frac{33}{32} > 1.$$
This completes step 3 of the proof.
\end{proof}

We return to our continuing example with $q=4$ and $n=3$. Recall that
$$F_{4,\le 3}(x) = 1 - \frac{x}3 + \frac{x^2}{45} - \frac{x^3}{2835},$$
and that its roots are
$$(r_{3,1},r_{3,2},r_{3,3}) \approx (3.997956, 16.80465, 42.19739).$$
The roots are alternately below and above $4$, $16$, and $64$; the first root
is very close to 4, the last one not so close.  We can expect this pattern
because, first, $$F_4(4^n) = 0$$ for any $n \ge 1$, and second, the difference
between $F_{4,\le 3}(x)$ and $F_4(x)$, 
$$F_{4,> 3}(x) = \frac{x^4}{722925} - \frac{x^5}{739552275} + \ldots,$$
is very small when $x$ is small.  The series $F_{4,> 3}(x)$ is also
dominated by its first term even when $x=64$.  Since $F_{4,> 3}(x)$ is positive
when $x \le 64$, the direction in which it displaces the first three roots of
$F_4(x)$ depends only on the sign of the derivative $F'_4(x)$. (The sign of the
derivative $f'(x)$ of any differentiable $f(x)$ must alternate between
consecutive simple roots.)  \lem{l:bound} is a careful estimate of the
displacement (indeed correct to within a universal constant factor).

To complete the proof of \thm{th:main}, we return to the definition of
$c_k$ in the statement of \lem{l:bound} and the convention $k=n+1-j$. Recall
that
$$a_{n,j} = \frac{1}{\sqrt{r_{n,j}}} \qquad q = p^2.$$
We claim that by \lem{l:bound},
\eq{e:final}{|p^{-j} - a_{n,j}| < p^{-j}\bigl|1-\frac{1}{\sqrt{c_k}}\bigr|
    < \frac2{p^{n+3k-1}} = \frac2{p^{j+4k-2}}.}
The first inequality is equivalent to \lem{l:bound}.   The second inequality is
far from sharp (\eqr{e:k2} is closer to the truth), but it is convenient to
prove \thm{th:main}.  To establish it, we need the elementary inequalities
$$1 - \frac{1}{\sqrt{1+x}} < \frac{x}2
    \qquad \frac{1}{\sqrt{1-x}}-1 < \frac{x}2 + \frac{3x^2}4$$
for $x>0$.  (They follow from the Taylor remainder theorem.)  When $k=1$, we
want to show that
$$\frac{1}{\sqrt{1-2p^{-2}}} - 1 < 2p^{-2}.$$
This can be established by symbolic algebra for $p \ge 2$ (indeed $p \ge
1.799$).  When $k>1$ is even,
$$1 - \frac{1}{\sqrt{1+4p^{-k(k+1)}}} < 2p^{-k(k+1)} \le 2p^{2-4k},$$
using that $k \ge 2$. When $k>1$ is odd,
\begin{align*}
\frac{1}{\sqrt{1-4p^{-k(k+1)}}} - 1 &< 2p^{-k(k+1)}+12p^{-2k(k+1)} \\
    &< 2p^{1-k(k+1)} < 2p^{2-4k},
\end{align*}
using that $p \ge 2$ and $k \ge 3$.

Finally the theorem follows from \eqr{e:final} by a
geometric sum:
$$(p-1)p^n\sum_{j=1}^n |a_{n,j}-p^{-j}| < \sum_{k=1}^n \frac{2(p-1)}{p^{3k-1}}
    < \frac{2p}{p^2+p+1} < 1,$$
as desired.

\section{Final remarks}

The proof of \lem{l:bound} obtains somewhat more information about the roots
$\{r_{n,j}\}$ of $F_{q,\le n}(x)$ than its statement. The proof shows that the
sequence
$$c_{n,k}= r_{n,j} q^{-j} = \frac{r_{n,n+1-k}}{q^{n+1-k}}$$
is monotonic in $n$ for every fixed $k$.  We can also change the bound $c_k$ to
be the solution to the equation
$$h_{q,k} = \frac{(-1)^{k+1}(1-c_k)q^{\binom{k+1}2}c_k^{-1-k}}
    {(c_k^{-1};q)_{-\infty}(c_k;q)_{-\infty}G_{q,\infty}(c_kq^{-k-1})} = 1$$
in the range $q^{-1/2} < c_k < q^{1/2}$.  (The equation is taken from
\eqr{e:hqk}.)
Then for this new value of $c_k$, 
$$\lim_{n \to \infty} c_{n,k} = c_k$$
and
$$\lim_{k \to \infty} (-1)^{k+1}(1-c_k)q^{\binom{k+1}2} = (1;q)_{-\infty}^3.$$
The value of this limit has an interesting interpretation when $q$ is a prime
power that may or may not be related to the present work. Its reciprocal is
the limiting probability that 3 independent, random $n \times n$ matrices over
the field $\F_q$ are non-singular.

That $c_{n,k}$ is monotonic in $n$ follows more directly from the interesting
recurrence
$$P_{q,\le n}(x) = (1-\frac{x}{q})P_{q,\le n-1}(\frac{x}{q})
    +\frac{(-x)^n}{q^n(q;q)_n}.$$
This recurrence also shows that $c_{n,k}$ is near $c_{n-1,k}$. This was the
basis of the author's first attempted proof of a lemma like \lem{l:bound}.
Such an attempt might yet have merit.

Finally we conjecture that \thm{th:main} together with its geometric
interpretation has a broad generalization to the mixed-base case:

\begin{conjecture} Let $p_1,p_2,\ldots,p_n \ge 2$ be a sequence of integers.
Let $X_1,X_2,\ldots,X_n$ be independent random variables such that $X_k$ is
uniformly distributed on the set
$$\{p_k-1,p_k-3,p_k-5,\ldots,1-p_k\}.$$
Then there are unique constants
$$a_1 > a_2 > \cdots > a_n > 0$$
such that the first $2n$ moments of
$$\tX = a_1 X_1 + a_2 X_2 + \cdots + a_n X_n$$
agree with the first $2n$ moments of a random variable $Y$ which is uniform on
$[-1,1]$. Moreover
$$\sum_{j=1}^n (p_j-1)\ \bigl|a_j- \prod_{k=1}^j p_k^{-1}\bigr|
    < \prod_{k=1}^n p_k^{-1}.$$
\eatline \label{c:main} \end{conjecture}

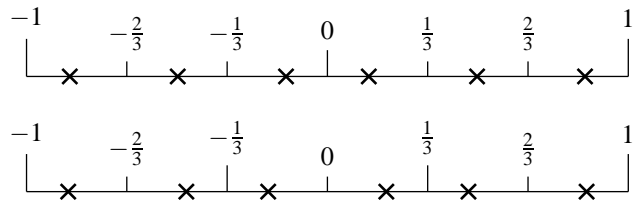
\begin{fullfigure}{f:ruler23}{The 6 values of $\tX$ marked on a ruler
    when $n=2$ and $(p_1,p_2) =(2,3)$ or $(p_1,p_2) =(3,2)$}
\pspicture(-4,0)(4,1.25)
\psline(-4,0)(4,0)
\psline(-4,0)(-4,.5)            \rput[b](-4,.65){\normalsize $-1$}
\psline(-2.667,0)(-2.667,.2)    \rput[b](-2.667,.35){\normalsize $-\frac23$}
\psline(-1.333,0)(-1.333,.2)    \rput[b](-1.333,.35){\normalsize $-\frac13$}
\psline(0,0)(0,.35)             \rput[b](0,.5){\normalsize $0$}
\psline(1.333,0)(1.333,.2)      \rput[b](1.333,.35){\normalsize $\frac13$}
\psline(2.667,0)(2.667,.2)      \rput[b](2.667,.35){\normalsize $\frac23$}
\psline(4,0)(4,.5)              \rput[b](4,.65){\normalsize $1$}
\rput(3.426608,0){\ptex} \rput(-3.426608,0){\ptex}
\rput(1.988709,0){\ptex} \rput(-1.988709,0){\ptex}
\rput(0.550810,0){\ptex} \rput(-0.550810,0){\ptex}
\endpspicture \\
\pspicture(-4,0)(4,1.5)  
\psline(-4,0)(4,0)       
\psline(-4,0)(-4,.5)            \rput[b](-4,.65){\normalsize $-1$}
\psline(-2.667,0)(-2.667,.2)    \rput[b](-2.667,.35){\normalsize $-\frac23$}
\psline(-1.333,0)(-1.333,.35)   \rput[b](-1.333,.5){\normalsize $-\frac13$}
\psline(0,0)(0,.2)              \rput[b](0,.35){\normalsize $0$}
\psline(1.333,0)(1.333,.35)     \rput[b](1.333,.5){\normalsize $\frac13$}
\psline(2.667,0)(2.667,.2)      \rput[b](2.667,.35){\normalsize $\frac23$}
\psline(4,0)(4,.5)              \rput[b](4,.65){\normalsize $1$}
\rput(3.445097,0){\ptex} \rput(-3.445097,0){\ptex}
\rput(1.874792,0){\ptex} \rput(-1.874792,0){\ptex}
\rput(0.785153,0){\ptex} \rput(-0.785153,0){\ptex}
\endpspicture
\eathalfline
\end{fullfigure}

For example, we can confirm Conjecture~\ref{c:main} when $n=2$ and either
$(p_1,p_2) = (2,3)$ or  $(p_1,p_2) = (3,2)$.  If we again let $r_j = a_j^{-2}$,
then in the first case,
$$r_1 = 15-2\sqrt{30} \qquad r_2 = 20+2\sqrt{30}$$
and
$$\tX \in \{\sim \pm .497177 \pm .179737,\sim \pm .497177\}.$$
In the second case,
$$r_1 = 20-2\sqrt{30} \qquad r_2 = 15+2\sqrt{30}$$
and
$$\tX \in \{\sim \pm .332493 \pm .196288,\sim \pm .196288\}.$$
The 6 values of $\tX$ in these two cases are
shown in \fig{f:ruler23}.

\acknowledgments

The author would like to thank Eric Rains for useful discussions. The author
would especially like to thank the referee for finding mistakes and unclear
points in the author's calculations. As the joke goes, the average paper is
read by one person, including both the author and the referee.  For the first
version of this paper, the referee was the one!


\providecommand{\bysame}{\leavevmode\hbox to3em{\hrulefill}\thinspace}

\end{document}